\documentclass[11pt, reqno]{amsart}
\usepackage{amsmath,amstext,amsthm,amsfonts}
\usepackage{amssymb,latexsym}
\usepackage[dvips]{graphicx}
\newtheorem{theorem}{Theorem }[section]

\theoremstyle{definition}
\theoremstyle{remark}

\newtheorem{claim}[theorem]{Claim}

\newcommand{\field}[1]{\mathbb{#1}}
\newcommand{\real}{\field{R}}

\newcommand{\de} {\delta}       \newcommand{\De}{\Delta}

\newcommand{\vep}{\varepsilon}

\newcommand{\te} {\theta}

\newcommand{\la} {\lambda}      \newcommand{\La}{\Lambda}

       \newcommand{\Si}{\Sigma}

       \newcommand{\Om}{\Omega}

\begin{document}
\title[LOCALLY CONFORMALLY FLAT MANIFOLDS]{CONVEXITY 
IN LOCALLY CONFORMALLY FLAT MANIFOLDS WITH BOUNDARY}

\author{Marcos Petr\'ucio de A. Cavalcante}


\noindent
\address{Instituto de Matem\'atica - Universidade Federal de Alagoas
\newline Campus A. C. Sim\~oes, BR 104 - Norte, Km 97, 57072-970.
\newline Macei\'o - AL -Brazil.
 }
\email{marcos@pos.mat.ufal.br}

\date{\today}


\begin{abstract}
{\bf Given a closed subset $\La$ of the open unit 
ball $B_1\subset \real^n$, $n \geq 3$, we will consider a complete 
Riemannian metric $g$ on $\overline{B_1} \setminus \La$ of 
constant scalar curvature equal to $n(n-1)$
and conformally related to the Euclidean metric. 
In this paper we prove that every closed Euclidean ball 
$\overline{B} \subset B_1\setminus \La$ is convex with respect to the metric $g$, 
assuming the mean curvature of the boundary $\partial B_1$ is nonnegative with respect to 
the inward normal. }
  \end{abstract}

\maketitle

\section{Introduction}
   Let $B_1$ denote the open unit ball of $\real^n$, $n \geq 3$. Given a
	closed subset $\La\subset B_1$, we will consider a complete Riemannian metric $g$ on
$\overline{B_1} \setminus \La$ of constant positive scalar curvature
$R(g)=n(n-1)$ and conformally related to the \index{Euclidean Metric, $\de$}Euclidean metric
$\de$. We will also assume that $g$ has nonnegative boundary mean curvature.
Here, and throughout this paper, second fundamental forms will be computed
with respect to the inward unit normal vector.

 In this paper we prove

\begin{theorem}\label{t1} 

If $B \subset B_1 \setminus \La$ is a standard Euclidean
ball, then $ \partial B$ is convex with respect to the metric $g$.
\end{theorem}

Here, we say that $\partial B$ is \index{Convex Boundary}\emph{convex} if its second
fundamental form is positive definite. Since $\partial B$ is
umbilical in the Euclidean metric and the notion of an umbilical
point is conformally invariant, we know that $\partial B$ is also
umbilic in the metric $g$. In that case $\partial B$ is convex if
its mean curvature $h$ is positive everywhere.

  This theorem is motivated by an analogous one on the sphere due
to R. Schoen ~\cite{S}. He shows that if $\La \subset S^n$ $n\geq3$, is
closed and nonempty and $g$ is a complete Riemannian metric on $S^n \setminus
\La$, conformal to the standard round metric $g_0$ and with constant
positive scalar curvature $n(n-1)$, then every standard ball $B \subset S^n
\setminus \La$ is convex with respect to the metric $g$. Schoen used
this geometrical result to prove the compactness of the set of solutions to
the Yamabe problem in the locally conformally flat case. Later, D.
Pollack also used Schoen's theorem to prove a compactness result for
the singular Yamabe problem on the sphere where the singular set is
a finite collection of points $\La = \{p_1,\dots,p_k\}\subset S^n$,
$n\geq 3$ (see ~\cite{P}).

In this context the Theorem \ref{t1} can be viewed as the 
first step in the direction of proving compactness for
the singular Yamabe problem with boundary conditions.

We shall point out that the problem of finding a metric satisfying 
the hypotheses of Theorem \ref{t1} is equivalent to finding a positive
solution to an elliptic PDE with critical Sobolev exponent. 
On the other hand this problem is invariant by conformal transformations.
So, by applying a convenient \emph{inversion} on the Euclidean space we 
may consider the same problem on an unbounded subset of $\real^n$.
The idea of the proof is to show that, if $\partial B$ is not convex,
then we can find a smaller ball $\widetilde B \subset B$ with non
convex boundary either. To do this we will use the hypothesis on the
mean curvature of $\partial B_1$ and get geometrical information from that 
equation by applying the Moving Planes Method as in ~\cite{GNN}.
The contradiction follows by the constructions of theses balls.



\section{Preliminaries}

  In this section we will introduce some notations and we shall recall
some results that will be used in the proof of Theorem \ref{t1}.
We will also describe a useful example.

Let $(M^n,g_0)$ be a smooth orientable Riemannian manifold, 
possibly with boundary, $n\geq3$.  Let us denote by  $R(g_0)$ 
its scalar curvature and by $h(g_0)$ its boundary mean curvature.
Let $g=u^{\frac{4}{n-2}}g_0$ be a \index{Conformal Metric} metric 
conformal to $g_0$. Then the positive function $u$ satisfies the 
following nonlinear elliptic partial differential equation of 
critical Sobolev exponent
\begin{equation}
 \label{edp}
   \left\{ \begin{array}{lrc}
  \De_{g_0} u - \frac{n-2}{4(n-1)}R(g_0) u +
  \frac{n-2}{4(n-1)}R(g)u^{\frac{n+2}{n-2}}=0 &  \textrm{ in } M, \\
  \frac{\partial u}{\partial \nu}-\frac{n-2}{2}h(g_0) u +
  \frac{n-2}{2}h(g)u^{\frac{n}{n-2}}=0 & \textrm{ on } \partial M,
   \end{array}  \right.
  \end{equation}
where $\nu$ is the inward unit normal vector field to $\partial M$.

The problem of existence of solutions to (\ref{edp}), when $R(g)$ and
$h(g)$ are constants, is referred to as the \index{Yamabe problem}\emph{Yamabe problem}.
It was completely solved when $\partial M = \emptyset$ in a sequence of
works, beginning with H. Yamabe himself ~\cite{Y}, followed by N. Trudinger
~\cite{T} and T. Aubin ~\cite{A}, and finally by R. Schoen ~\cite{S2}.
In the case of nonempty  boundary, J. Escobar solved almost all the cases
(see ~\cite{E1}, ~\cite{E2}) followed by Z. Han and Y. Li ~\cite{HL},
F. Marques ~\cite{M} and others. In this article, however, we wish to
study solutions of (\ref{edp}), with $R(g)$ constant, which become
singular on a closed subset $\La \subset M$. This is the so
called \index{Singular Yamabe Problem}\emph{singular Yamabe problem}.
This singular behavior is equivalent,
at least in the case that $g_0$ is conformally flat, to requiring
$g$ to be complete on $M \setminus \La$.  The existence problem (with
$\partial M=\emptyset$) displays a relationship between the size of
$\La$ and the sign of $R(g)$. It is known that for a solution with
$R(g)<0$ to exist, it is necessary and sufficient that $\textrm{dim}
(\La)> \frac{n-2}{2}$ (see ~\cite{AMc}, ~\cite{Mc} and ~\cite{FMc}),
while  if a solution exists with $R(g) \geq 0$, then $\textrm{dim}(\La)
\leq \frac{n-2}{2}$. Here $\textrm{dim}(\La)$ stands for the Hausdorff
dimension of $\La$. In this paper we will treat the case of constant
positive scalar curvature, which we suppose equal to $n(n-1)$ after 
normalization. In this case the simplest examples are given by the Fowler
solutions which we will now discuss briefly.

\vspace{0.2cm}

Let $u:\real^n\setminus \{0\}\to\real$ be a positive
smooth function such that 

\begin{equation}
 \label{eq_rn}
   \left\{ \begin{array}{lcc}
  \De u +  \frac{n(n-2)}{4}u^{\frac{n+2}{n-2}}=0 &  \textrm{ in }
  \real^n \setminus\{0\}, n\geq 3, \\
   0 \textrm{ is an isolated singularity}.
   \end{array}  \right.
  \end{equation}
 \index{Conformal Metric!Equation in $\real^n$}

In this case, $g=u^{\frac{4}{n-2}}\de$ is a complete metric on 
$\real^n \setminus \{0\}$ of constant scalar
curvature ~$n(n-1)$.
%




%

Using the invariance under conformal transformations we may work in
different background metrics. The most convenient one here
is the \index{Cylindrical Metric}cylindrical metric $g_{cyl}=d\te^2+dt^2$ on $S^{n-1} \times \real$.
In this case $g=v^{\frac{4}{n-2}}g_{cyl}$, where $v$  is defined in the
whole cylinder and satisfies

\begin{equation}
 \label{edp_cyl}
  \frac{d^2v}{dt^2}+ \De_{\te} v -\frac{(n-2)^2}{4}v +
\frac{n(n-2)}{4}v^{\frac{n+2}{n-2}}=0.
  \end{equation}

One easily verifies that  the solutions  to equation (\ref{eq_rn}) and
(\ref{edp_cyl}) are related by

\begin{equation}
\label{relation}
u(x)=|x|^{\frac{2-n}{2}}v(x/|x|,-\log |x|).
\end{equation}

By a deep theorem of Caffarelli, Gidas and Spruck (see ~\cite{CGS}, Theorem 8.1)
we know that $v$ is rotationally symmetric, that is $v(\te,t)= v(t)$, and therefore the
PDE (\ref{edp_cyl}) reduces to the following ODE:

\begin{eqnarray*}
  \frac{d^2v}{dt^2} -\frac{(n-2)^2}{4}v +
\frac{n(n-2)}{4}v^{\frac{n+2}{n-2}}=0.  \nonumber
  \end{eqnarray*}

Setting $w=v'$ this equation is transformed into a first order
\index{Hamiltonian System}Hamiltonian system

\begin{eqnarray*}
   \left\{ \begin{array}{l}
 \frac{dv}{dt}= w,\\
 \frac{dw}{dt}= \frac{(n-2)^2}{4}v -
\frac{n(n-2)}{4}v^{\frac{n+2}{n-2}},
  \nonumber \\
      \end{array}   \right.
  \end{eqnarray*}
whose Hamiltonian energy is given by

\begin{eqnarray*}
H(v,w)= w^2 -\frac{(n-2)^2}{4}v^2+ \frac{(n-2)^2}{4}v^{\frac{2n}{n-2}}.
  \nonumber 
\end{eqnarray*}

The solutions $(v(t),v'(t))$ describe the level sets of $H$ and we
note that $(0,0)$ and $(\pm v_0,0)$, where $v_0= \big(\frac{n-2}{n}
\big)^{\frac{n-2}{4}}$, are the equilibrium points. We restrict ourselves
to the half-plane  $\{v>0\}$ where $g=v^{\frac{4}{n-2}}g_{cyl}$ has
geometrical meaning. On the other hand we are looking for complete
metrics. Those will be generated by the \index{Fowler Solutions}\emph{Fowler solutions}: the
periodic solutions around the equilibrium point $(v_0,0)$. They are
symmetric with respect to $v$-axis and
can be parametrized by the minimum value $\vep$ attained by $v$,
$\vep \in (0, v_0]$, (and a translation parameter $T$).
We will denote them by $v_{\vep}$. We point out
that $v_0$ corresponds to the scaling of $g_{cyl}$ which makes
the cylinder $S^{n-1}\times \real$ 
have scalar curvature $n(n-1)$. We observe that one obtains
the Fowler solutions $u_\vep$ in $\real^n\setminus\{0\}$ by using the relation
(\ref{relation}).

We can now construct metrics satisfying the hypotheses of Theorem
~\ref{t1} (with $\La =\{0\}$) from the Fowler solutions.
To do this, we just take a Fowler solution $v$ defined
for $t \geq t_0$, where $t_0$ is such that we have $w=\frac{dv}{dt}
 \leq 0$, or equivalently,
\begin{eqnarray*}
h(g)= -\frac{2}{n-2}v^{-\frac{n}{n-2}} \frac{dv}{dt} \geq 0.
\end{eqnarray*}
\vspace{.2cm}

We point out that, by another result of Caffarelli, Gidas and Spruck
(see Theorem 1.2 in ~\cite{CGS}) it is known that, given a positive solution
$u$ to
\begin{eqnarray}
\label{eq}
  \De u +  \frac{n(n-2)}{4}u^{\frac{n+2}{n-2}}=0
\end{eqnarray}
which is defined in the punctured ball $B_1\setminus \{0\}$ and which
is singular at the origin, there exists a unique Fowler solution $u_\vep$
such that
$$
u(x)=(1+o(1))u_\vep(|x|) \, \textrm{ as } \, |x|\to 0.
$$
Therefore, from equation (\ref{relation}) (see also ~\cite{KMPS}), 
either $u$ extends as a smooth solution to the ball, or there exist 
positive constants $C_1$, $C_2$ such that
\begin{eqnarray*}
\label{assint}
C_1|x|^{(2-n)/2} \leq u(x)\leq C_2|x|^{(2-n)/2}.
\end{eqnarray*}

\section{Proof of Theorem \ref{t1} }

The proof will be by contradiction. 
If $\partial B$ is not convex then, since it is umbilical, 
there exists a point $q \in \partial B$ such that the mean curvature
of $\partial B$ at $q$ (with respect to the inward unit normal vector) is $H(q)\leq 0$.
If we write $g=u^{\frac{4}{n-2}}\de$ we have that $u$ is a positive
smooth function on $\overline{B}_1 \setminus \La$ satisfying
 \begin{equation}
   \left\{ \begin{array}{lr}
   \label{y_eq_in_b}
  \De u +  \frac{n(n-2)}{4}u^{\frac{n+2}{n-2}}=0 &  \textrm{ in } B_1
\setminus \La, \\
  \frac{\partial u}{\partial \nu} - \frac{n-2}{2}u +
    \frac{n-2}{2}h u^{\frac{n}{n-2}}=0 & \textrm{ on } \partial B_1. \\
 \end{array}    \right.
\end{equation}


Now, we will choose a point $p\in \partial B$, $p \neq q$ and let us
consider the inversion
\begin{eqnarray*}
I:\real^n\setminus \{p\}\rightarrow \real^n\setminus \{p\}.
\end{eqnarray*}

This map  takes $\overline {B_1} \setminus (\{p\}\cup\La$)
on $\real^n \setminus (B(a,r) \cup \La)$, where $B(a,r)$
is an open ball of center $a \in \real ^n$ and radius $r>0$ and
$\La$ still denotes the singular set. Let us denote by $\Si$ the
boundary of $B(a,r)$, that is, $\Si = I(\partial B_1)$.

The image of $\partial B\setminus \{p\}$ is a hyperplane $\Pi$ and by a coordinate
choice we may assume $\Pi=\Pi_0:=\{x\in \real^n :  x^n=0\}$. We may
suppose that the ball $B(a,r)$ lies below $\Pi_0$. 
Notice that in this case $\La$ also lies below $\Pi_0$.

Since $I$ is a conformal map we have $I^*g= v^{\frac{4}{n-2}}\de$,
where $v$ is the \emph{Kelvin transform} of $u$ on $\real^n \setminus
(B(a,r) \cup \La)$.

Thus this metric has constant positive scalar curvature $n(n-1)$ in
$\real^n \setminus (B(a,r) \cup \La)$
and nonnegative mean curvature $h$ on $\Si$.

As before $v$ is a solution of the following problem
\begin{eqnarray*}
  \left\{ \begin{array}{lr} \De v +
\frac{n(n-2)}{4}v^{\frac{n+2}{n-2}}=0 &
   \textrm{ in  } \real^n \setminus (B(a,r)
  \cup \La),\\
  \frac{\partial v}{\partial
  \nu}+\frac{n-2}{2r}v+\frac{n-2}{2}hv^{\frac{n}{n-2}}=0 & \textrm{
on    } \Si.
  \end{array}   \right.
\end{eqnarray*}

Also, by hypotheses of contracdition, the mean curvature of
the hyperplane $\Pi_0$ at $I(q)$ (with respect to $\frac{\partial}{\partial x^n}$ )
is $H\leq 0$. By applying the boundary equation of the system (\ref{edp}) 
to $\Pi_0$ we obtain $\frac{\partial v}{\partial x^n} + \frac{n-2}{2}H
v^{\frac{n}{n-2}}=0$ on $\Pi_0$. Thus we conclude that 
$\frac{\partial v}{\partial x^n} (I(q))\geq 0$.



\vspace{0.15cm}
Now we start with the Moving Planes Method.
Given $\la \geq 0$ we will denote by $x_{\la}$
the reflection of $x$ with respect to the hyperplane
$\Pi_{\la}:=\{x\in \real^n: x^n=\la\}$ and set $\Om_\la= \{x \in \real^n
\setminus(B(a,r) \cup \La): x^n\leq \la\}$. We define
$$
w_\la(x)=v(x)-v_\la(x)\textrm{ for }x\in \Om_\la,
$$
where $v_\la(x):=v(x_\la)$.

Since the infinity is a regular point of $I^*g$, we have that
$$
v(x)=|x|^{2-n}\Big(a+ \sum b_ix^i|x|^{-2}\Big) +O(|x|^{-n})
$$
in a neighborhood of infinity.
It follows from Lemma 2.3 of \cite{CGS} 
that there exist $R>0$ and $\bar{\la}>0$ such that $w_\la>0$ in
interior of 
$\Om_\la \setminus B(0,R)$, if $\la\geq\bar{\la}$. Without loss of
generality we can choose $R>0$ such that $B(a,r) \cup \La
\subset B(0,R)$.

Now we note that $v$ has a positive infimum, say $v_0>0$, in  $B(0,R) \setminus
(B(a,r)\cup\La)$. It follows from the fact that $v$ is a classical
solution to (\ref{eq}) in $B(0,R) \setminus (B(a,r)\cup\La)$.
So, since $v$ decays in a neighborhood of infinity, 
we may choose $\overline{\la}>0$
large enough such that $v_\la(x)<v_0/2$,
for $x \in B(0,R)$ and for $\la \geq \overline{\la}$.
Thus, for sufficiently large $\la$ we get $w_\la > 0$ in
$\textrm{int}(\Om_\la)$. 

We also write
  \begin{equation}
  \label{eq_for_w_lamba}
  \De w_\la + c_{\la}(x)w_\la=0 \textrm{ in } \textrm{int}(\Om_\la),
  \end{equation}
where
$$
c_{\la}(x)=\frac{n(n-2)}{4}\frac{v(x)^{\frac{n+2}{n-2}}-
v_\la(x)^{\frac{n+2}{n-2}}}{v(x)-v_\la(x)}.
$$

Notice that, by definition, $w_\la$ always vanishes on $\Pi_\la$. In
particular, setting $\la_0=\inf\{\overline{\la}>0: w_\la > 0\textrm{ on }
\textrm{int}(\Om_\la), \forall \la\geq \overline{\la} \}$ we obtain by
continuity that  $w_{\la_0}$ satisfies (\ref{eq_for_w_lamba}),
$w_{\la_0}\geq 0$ in $\Om_{\la_0}$ and $w_{\la_0}=0$ on $\Pi_{\la_0}$.
Hence, by applying the strong maximum principle, we conclude that either
$w_{\la_0}>0$ in $\textrm{int}(\Om_{\la_0})$ or $w_{\la_0}=v-v_{\la_0}$ vanishes
identically. We point out that the second case occurs only if $\La = \emptyset$.

If $w_{\la_0} \equiv 0$, then $\Pi_{\la_0}$ is a hyperplane of symmetry 
of $v$ and therefore $v$ extends to a global positive 
solution of (\ref{eq}) on the entire $\real^n$. Using ~\cite{CGS}, we 
conclude that $(B_1, g)$ is a convex spherical cap and the result is obvious.

If $w_{\la_0}>0$ in $\textrm{int}(\Om_{\la_0})$ we apply the E. Hopf maximum principle
to conclude
  \begin{eqnarray}
  \label{hopflemma}
  \frac{\partial w_{\la_0}}{\partial x^n}=2 \frac{\partial v}{\partial
x^n}<0
  \textrm{ in } \Pi_{\la_0},
  \end{eqnarray}
and since $\frac{\partial v}{\partial x^n}(I(q))\geq 0$, we have
$\la_0 >0$. In this case, by definition of $\la_0$, we can choose
sequences $\la_k \uparrow \la_0$ and
$x_k \in \Om_{\la_k}$ such that $w_{\la_k}(x_k)< 0$.

It follows from the work in \cite{KMPS} that
$w_\la$ achieves its infimum.
%
%
Then we may assume, without loss of generality,
that $x_k$ is a minimum of $w_{\la_k}$ in $\Om_{\la_k}$.

We have that $x_k \notin \Pi_k$ because $w_{\la_k}$ always vanishes on $\Pi_{\la_k}$.
So, either $x_k\in\Si$ or is an interior point.
Even when $x_k$ is an interior point we claim that $(x_k)_k$ is a bounded
sequence. More precisely,

\begin{claim}
\label{claimappendix}[see \S 2 in ~\cite{CL}] There exists $R_0>0$,
independent of $\la$, such that if $w_{\la}$ solves (\ref{eq_for_w_lamba}) and 
is negative somewhere in 
$\textrm{int}(\Om)$, and $x_0\in
\textrm{int}(\Om)$ is a minimum point of $w_{\la}$, then $|x_0|<R_0$.
\end{claim}

For completeness we present a proof in the Appendix.

So, we can take a convergent subsequence  $x_k \to\bar{x}\in\Om_{\la_0}$.
Since $w_{\la_k}(x_k)<0$ and $w_{\la_0}\geq 0$ in $\Om_{\la_0}$ we necessarily
have $w_{\la_0}(\bar{x})=0$ and therefore 
$\bar{x} \in \partial \Om_{\la_0}= \Pi_{\la_0}\cup \Si$. 

If $x\in \Pi_{\la_0}$ then $x_k$ is an interior minimum point to $w_{\la_k}$ 
and hence $\nabla w_{\la_0}(\bar{x})=0$ which not ocurrs by 
inequality (\ref{hopflemma}). Thus 
we have $\bar{x} \in \Si$ and by E. Hopf maximum principle again,
\begin{equation}
\label{derivadadew}
\frac{\partial w_{\la_0}}{\partial \eta}(\bar{x})=
\frac{\partial v}{\partial \eta}(\bar{x}) -
\frac{\partial v}{\partial \eta}(\bar{x}_{\la_0}) <0,
\end{equation}
where $\eta:=-\nu$ is the inward unit normal vector to $\Si$.

%
%

Now, we recall that
\begin{equation} \label{eqdevnobord}
\frac{\partial v}{\partial \nu} + \frac{n-2}{2r}v+\frac{n-2}{2}hv^{\frac{n+2}{n-2}}=0
 \,\textrm{ on } \, \Si.
\end{equation}

Thus, since $v(\bar{x})=v(\bar{x}_{\la_0})$ we have from (\ref{derivadadew}) and
(\ref{eqdevnobord}) that the mean curvature of $\Si_{\la_0}$ at $\overline{x}_{\la_0}$
(with respect to the inward unit normal vector) is strictly less than $-h$.

Since $h\geq 0$, we have that
$\bar{x}_{\la_0}$ is a non convex point in the reflected sphere $\Si_{\la_0}$
Considering the problem
back to $B_1$ we denote by $K_1$ the ball corresponding to the ball whose boundary is
$\Si_{\la_0}$ and by $P_1$ the ball corresponding to $\Pi_{\la_0}^+ $. 
Thus we have obtained a strictly smaller ball $K_1 \subset B$ with 
non convex boundary which is the reflection of $\partial B_1 $
with respect to $\partial  P_1$.

We can repeat this argument to obtain a sequence of balls with non
convex points on the boundaries,
$B  \supset K_1 \supset \cdots  \supset K_j \supset \cdots$.

This sequence cannot converge to a point, since
small balls are always convex. On the other hand, if  $K_j
\rightarrow K_{\infty}$ where $K_{\infty}$ is not a point, then
$K_{\infty} \subset B$ is a  ball in $B_1 \setminus \La$
such that its boundary is the reflection of
$\partial B_1$  with respect to to itself, that is a contradiction.

%


\appendix

\section{Proof of Claim \ref{claimappendix}}

First write (\ref{eq_for_w_lamba}) setting $c_{\la}(x) = 0$ when
$w_{\la}(x)=0$.  Fix  $0<\mu < n-2$ and define $g(x)=|x|^{-\mu}$ and
$\phi (x)=\frac{w_\la (x)}{g(x)}$. Then, using the equation
(\ref{eq_for_w_lamba}),

$$
\De \phi +\frac{2}{g}\langle \nabla g,\nabla \phi \rangle + \left(
c_{\la}(x)+ \frac{\De g}{g}\right)\phi=0.
$$
By a computation we get $\De g = -\mu(n-2-\mu)|x|^{-\mu -2}$, that
is,

  $$\frac{\De g}{g}= -\mu(n-2-\mu)|x|^{-2}.$$

  On the other hand, the expansion of $v$ in a neighborhood of infinity
implies that $w_{\la}(x)= O(|x|^{2-n})$ and
consequently  $c_{\la}(x) = O(|x|^{-n-2-2+n})= O(|x|^{-4})$.  Hence we
obtain

$$
c_{\la}(x)+ \frac{\De g}{g} \leq C(|x|^{-4}- \mu (n-2-\mu))|x|^{-2}).
$$

In particular $c(x)+ \frac{\De g}{g} < 0$ for large $|x|$. Choose
$R_0$ with  $B(a, r)\cup \La \subset B(0, R_0)$ such that
\begin{equation}
\label{x_less_zero}
  C(|x|^{-4}- \mu (n-2-\mu))|x|^{-2})< 0, \textrm{ for } |x| \geq R_0.
\end{equation}

 Now let $x_0 \in \textrm{int}(\Om_\la)$ so that $w_\la(x_0)=
\inf_{\textrm{int}(\Om_\la)}w_\la <0$.

Since $\lim _{|x| \rightarrow +\infty}\phi(x)=0$ and $\phi(x) \geq 0$
on $\partial \Om_\la$, there exists $\bar{x}_0$ such that $\phi$ has its
minimum at $\bar{x}_0$. By applying the maximum principle for $\phi$ at
$\bar{x}_0$ we get $c_{\la}(\bar{x}_0) + \frac{\De g(\bar{x}_0)}{g} \geq 0$
and by (\ref{x_less_zero}), $|\bar{x}_0| < R_0$. Now we have
\begin{eqnarray*}
  \begin{array}{ccc}
\frac{w_{\la}(x_0)}{g(\bar{x}_0)} &\leq&
\frac{w_{\la}(\bar{x}_0)}{g(\bar{x}_0)} =\phi (\bar{x}_0)\\
&\leq& \phi (x_0) =\frac{w_\la(x_0)}{g(x_0)}.
   \end{array}
\end{eqnarray*}
This implies $|x_0|\leq |\bar{x}_0| \leq R_0$ and proves the claim.

\vspace{0.1cm}

\noindent
{\bf Acknowledgements:} 
The content of this paper is part of the author's doctoral thesis \cite{C}. 
The author would like to express his gratitude to Prof. Manfredo do Carmo 
for the encouragement and to  Prof. Fernando Coda Marques for many 
useful discussions during this work. While the author was at IMPA - Rio de Janeiro, 
he was fully support by CNPq-Brazil.

\end{document}